\documentclass[11pt]{article}

\usepackage[T1]{fontenc}
\usepackage{lmodern}

\usepackage{amsmath, amssymb, amsthm}
\usepackage{mathrsfs}
\usepackage{stmaryrd}
\usepackage[margin=2.5cm]{geometry}
\usepackage{color, graphicx, float}
\usepackage[breakable]{tcolorbox}

\usepackage{etoolbox}
\makeatletter
\patchcmd{\@maketitle}{\LARGE \@title}{\LARGE\bfseries\@title}{}{}

\renewcommand{\@seccntformat}[1]{\csname the#1\endcsname.\quad}
\makeatother

\definecolor{darkblue}{rgb}{0,0,.5}

\usepackage{hyperref}
\hypersetup{
	colorlinks=true,		
	linkcolor=darkblue,		
	citecolor=darkblue,		
	urlcolor=darkblue 		
}

\makeatletter
\def\th@plain{%
	\thm@notefont{}
	\itshape 
}
\def\th@definition{%
	\thm@notefont{}
	\normalfont 
}

\renewenvironment{proof}[1][\proofname]{\par
	\normalfont
	\topsep0\p@\@plus3\p@ \trivlist
	\item[\hskip\labelsep\itshape
	#1\@addpunct{.}]\ignorespaces
}{%
	\qed\endtrivlist
}
\renewcommand\paragraph{\@startsection{paragraph}{4}{\z@}%
    {1.5ex \@plus1ex \@minus.2ex}%
    {-1em}%
    {\normalfont\normalsize\bfseries}}
\makeatother

\newtheorem{theorem}{Theorem}[section]
\newtheorem{lemma}[theorem]{Lemma}

\theoremstyle{definition}
\newtheorem{definition}[theorem]{Definition}
\theoremstyle{definition}
\newtheorem{example}[theorem]{Example}
\theoremstyle{definition}
\newtheorem{remark}[theorem]{Remark}
\theoremstyle{definition}
\newtheorem{algorithm}{Algorithm}

\usepackage[shortlabels]{enumitem}

\allowdisplaybreaks

\newcommand{\argmin}{\ensuremath{\operatorname*{argmin}}}
\newcommand{\dom}{\ensuremath{\operatorname{dom}}}

\newcommand{\epi}{\ensuremath{\operatorname{epi}}}
\newcommand{\Id}{\ensuremath{\operatorname{Id}}}
\newcommand{\prox}{\ensuremath{\operatorname{Prox}}}
\newcommand{\Limsup}{\ensuremath{\operatorname*{Limsup}}}

\newcounter{step}[algorithm]
\newcommand\step[1]{%
	\refstepcounter{step}	
	\vskip 0.25\baselineskip
	\ifx\hfuzz#1\hfuzz
		\item[~\(\triangleright\)~\textbf{Step~\arabic{step}.}]
	\else
		\item[~\(\triangleright\)~\textbf{Step~\arabic{step}}] (\texttt{#1})\textbf{.}%
	\fi
}

\begin{document}

\title{A proximal subgradient algorithm for constrained multiobjective DC-type optimization}

\author{
Nguyen Van Tuyen\thanks{Department of Mathematics, Hanoi Pedagogical University 2, Xuan Hoa, Phu Tho, Vietnam.
E-mail:~\href{href:nguyenvantuyen83@hpu2.edu.vn; tuyensp2@yahoo.com}{nguyenvantuyen83@hpu2.edu.vn; tuyensp2@yahoo.com}.},
~
Minh N. Dao\thanks{School of Science, RMIT University, Melbourne, VIC 3000, Australia.
E-mail:~\href{href:minh.dao@rmit.edu.au}{minh.dao@rmit.edu.au}.},
~and~
Tran Van Nghi\thanks{Department of Mathematics, Hanoi Pedagogical University 2, Xuan Hoa, Phu Tho, Vietnam.
E-mail:~\href{href:tranvannghi@hpu2.edu.vn}{tranvannghi@hpu2.edu.vn}.}
}

\date{September 3, 2026}

\maketitle

\begin{abstract}
In this paper, we consider a class of constrained multiobjective optimization problems, where each objective function can be expressed by adding a possibly nonsmooth nonconvex function and a differentiable function with Lipschitz continuous gradient, then subtracting a weakly convex function. This encompasses multiobjective optimization problems involving difference-of-convex (DC) functions, which are prevalent in various applications due to their ability to model nonconvex problems. We first establish a necessary optimality condition for these problems and then derive sufficient optimality conditions under some structural assumptions, providing a theoretical foundation for algorithm development. Building on these conditions, we propose a proximal subgradient algorithm tailored to the structure of the objectives. Under mild assumptions, the sequence generated by the proposed algorithm is bounded and each of its cluster points is a stationary point. Numerical experiments on multi-sensor sparse recovery illustrate the effectiveness and computational advantages of the proposed algorithm.  
\end{abstract}

\noindent{\bfseries Keywords:}
DC programming,
multiobjective optimization,
multi-sensor sparse recovery,
optimality conditions,
proximal subgradient algorithm.

\noindent{\bf Mathematics Subject Classification (MSC 2020):}
90C26,  
90C29,  
90C46,  
65K05.  

\section{Introduction}

Multiobjective optimization is a field that seeks to optimize multiple conflicting objectives simultaneously, often requiring trade-offs among the objectives due to the absence of a universally optimal solution. Such problems are fundamental in diverse domains, including economics, engineering, and the sciences, where decision-making processes require balancing competing goals effectively; see, e.g., \cite{Ehrgott,Jahn-04,Luc-89,Sawaragi-et al-85}.

A variety of approaches have been developed for solving multiobjective optimization problems. Classical scalarization techniques transform the vector problem into one or a family of scalar problems, for instance through weighted sums or minimax scalarizations; see \cite{Ehrgott}. Although flexible and widely used, such approaches often require the selection of scalarization parameters and may target only particular trade-offs among the objectives. Alternatively, descent methods seek directions that simultaneously improve the objective components without prescribing weights in advance; a representative example is the steepest descent method in \cite{Fliege-Svaiter-00}. Proximal regularization provides another important methodology for vector objectives; see, e.g., \cite{Bonnel-Iusem-Svaiter-05}. For composite multiobjective optimization, proximal gradient methods, with and without line searches, are developed in \cite{Tanabe-Fukuda-Yamashita-19} for problems in which each objective is the sum of a continuously differentiable function and a proper convex nonsmooth function.

Difference-of-convex (DC) structure provides a useful framework for modeling nonsmooth and nonconvex optimization problems and has been extensively studied in scalar optimization. Its broad modeling capacity and amenability to convex-analytic techniques have led to applications in areas such as machine learning and engineering \cite{An-Tao-18,de Oliveira}. Theoretical aspects of DC optimization, including optimality conditions, duality, and polyhedral DC structures, have also been extensively investigated; see, e.g., \cite{Bomze-10,Hang-Yen-16,Huong-Huyen-Yen-24,Tao-An}. In the multiobjective setting, proximal point methods for problems with DC objective functions are developed in \cite{Ji-Goh-Souza-16}, while a proximal approach for constrained multiobjective DC optimization is proposed in \cite{Bento-et-al-20}. These developments show that proximal regularization can be effectively combined with DC decompositions in multiobjective optimization. However, such approaches are primarily designed for problems in which each objective is represented as the difference of two convex functions and therefore do not directly exploit situations in which the individual components of an objective possess different smoothness and convexity properties.

In this paper, we consider a class of constrained multiobjective optimization problems where each objective function is formed by adding a possibly nonsmooth nonconvex locally Lipschitz function and a differentiable function with Lipschitz continuous gradient, then subtracting a weakly convex function. This general formulation encompasses important classes of multiobjective optimization problems, including those with DC objective functions and those with additive composite structures. At the same time, it allows the individual components of each objective to have substantially different analytical properties. This broader structure creates several challenges. In particular, the presence of a nonconvex locally Lipschitz term and a possibly nonsmooth weakly convex subtractive term requires optimality analysis adapted to the more general setting. From an algorithmic perspective, the three distinct nonsmooth, smooth, and weakly convex components should also be treated in a way that exploits their respective properties. The presence of a constraint set further complicates both the optimality analysis and the construction of computationally tractable descent steps.

Our contributions address both theoretical and computational aspects of the considered problem class. We first establish a necessary optimality condition for local weak Pareto solutions and derive sufficient optimality conditions under additional structural assumptions. These results provide optimality characterizations for the considered constrained multiobjective setting, including cases where the subtractive functions are nonsmooth. A proximal subgradient algorithm is then proposed to exploit the three-part composite structure of the objectives through a max-type model of their centered local variations, without fixing scalarization weights for the original vector objective in advance. We demonstrate a simultaneous descent estimate for all objective components and, under suitable assumptions, boundedness of the generated sequence and stationarity of its cluster points.

The remainder of the paper is organized as follows. Section~\ref{s:prelim} provides the preliminaries, including key definitions and concepts relevant to multiobjective DC-type optimization. In Section~\ref{s:optcond}, we derive a necessary optimality condition and sufficient optimality conditions under structural assumptions. Section~\ref{s:algo} introduces the proposed proximal subgradient algorithm and presents its convergence analysis. In Section~\ref{s:numerical}, we report numerical experiments on multi-sensor sparse recovery to illustrate the practical performance of the proposed algorithm. Finally, Section~\ref{s:conclu} concludes the paper and discusses possible directions for future research.

\section{Preliminaries}
\label{s:prelim}

Let us recall some notions related to generalized differentiation from~\cite{Mordukhovich2006,Mordukhovich2018,Rockafellar1998}. Throughout the paper, we work in the Euclidean space $\mathbb{R}^n$ equipped with the usual inner product $\langle \cdot, \cdot\rangle$ and the corresponding norm $\|\cdot\|$. 

For a set-valued mapping $F\colon \mathbb{R}^n\rightrightarrows \mathbb{R}^m$, the limiting construction
\begin{align*}
\Limsup_{x\to \bar x} F(x):=\bigg\{ y\in \mathbb{R}^m \mid \exists x_k \to \bar x, y_k \to y \ \mbox{with}\ y_k\in F(x_k)\bigg\}
\end{align*}
is known as the \emph{Painlev\'e--Kuratowski outer/upper limit} of $F$ at $\bar x$.

\begin{definition}[See \cite{Mordukhovich2006,Mordukhovich2018}]
Let $S$ be a nonempty subset of $\mathbb{R}^n$ and $\bar x \in S$.
\begin{enumerate}
\item
The \emph{regular/Fr\'echet normal cone} to $S$ at $\bar x$ is defined by
\begin{align*}
\widehat N(\bar x, S) :=\left\{ v\in \mathbb{R}^n\mid \limsup\limits_{x \xrightarrow{S}\bar x} \dfrac{\langle v, x-\bar x \rangle}{\|x-\bar x\|} \leq 0 \right\},
\end{align*}
where $x \xrightarrow{S} \bar x$ means that $x \to \bar x$ and $ x\in S$.
\item
The \emph{limiting/Mordukhovich  normal cone} to $S$ at $\bar x$ is given by
\begin{align*}
N(\bar x, S) :=\Limsup\limits_{ x \xrightarrow{S} \bar x} \widehat{N}(x, S).
\end{align*}
\end{enumerate}
We set $\widehat N(\bar x,S) =N(\bar x,S) :=\varnothing$ if $\bar x \not\in S$.
\end{definition}

Given a function $\varphi\colon \mathbb{R}^n\to (-\infty, +\infty]$, its \emph{effective domain} is
\begin{align*}
\dom \varphi:=\{x \in \mathbb{R}^n \mid \varphi(x) < +\infty\}    
\end{align*}
and its \emph{epigraph} is
\begin{align*}
\epi \varphi:=\{ (x, \alpha) \in \mathbb{R}^n \times \mathbb{R} \mid \alpha \ge \varphi (x)\}.    
\end{align*}
We call $\varphi$ a \emph{proper} function if $\dom \varphi$ is a nonempty set. The function $\varphi$ is said to be \emph{locally Lipschitz around a given point $\bar x\in\mathbb{R}^n$} if there is a neighborhood $U\subseteq \dom \varphi$ of $\bar x$ and a constant $\ell\geq 0$ such that, for all $x, y\in U$,
\begin{align*}
|\varphi(x)-\varphi(y)|\leq \ell \|x-y\|.
\end{align*}
We say that $\varphi$ is \emph{locally Lipschitz on an open subset $G$} of $\mathbb{R}^n$ if it is locally Lipschitz around every point of this set.

\begin{definition}[See \cite{Mordukhovich2006,Mordukhovich2018}]
Consider a function $\varphi \colon \mathbb{R}^n \to (-\infty, +\infty]$ and a point $\bar{x} \in \dom \varphi$.
\begin{enumerate}
\item
The \emph{regular/Fr\'echet subdifferential} of $\varphi$ at $\bar{x}$ is
\begin{align*}
\widehat{\partial}\varphi(\bar{x}):=\{ v \in \mathbb{R}^n \mid (v,-1)\in \widehat{N}_{\epi \varphi}(\bar{x},\varphi(\bar{x})) \}.    
\end{align*}
\item
The \emph{limiting/Mordukhovich subdifferential} of $\varphi$ at $\bar{x}$ is
\begin{align*}
\partial \varphi(\bar{x}):=\{ v \in \mathbb{R}^n \mid (v,-1)\in {N}_{\epi \varphi}(\bar{x},\varphi(\bar{x})) \}.    
\end{align*}
\end{enumerate}
\end{definition}

It is well known that
\begin{align*}
\partial \varphi(\bar{x})=\Limsup_{x \xrightarrow{\varphi} \bar{x}}\widehat{\partial}\varphi(x) \supseteq \widehat{\partial}\varphi(x),
\end{align*}
where $x \xrightarrow{\varphi} \bar{x}$ means that $x \to \bar{x}$ and $\varphi(x) \to \varphi(\bar{x})$. In
particular, if $\varphi$ is a convex function, then the subdifferentials $\widehat{\partial} \varphi(\bar{x})$ and $\partial \varphi(\bar{x})$
coincide with the subdifferential in the sense of convex analysis.

Let $S \subseteq \mathbb{R}^n$. We define the \emph{indicator function} $\delta_{S} \colon \mathbb{R}^n \to (-\infty, +\infty]$ by
\begin{align*}
\delta_{S}(x) :=
\begin{cases}
0 & \textrm{ if } x \in S, \\
+\infty & \textrm{ otherwise.}
\end{cases}
\end{align*}
It holds that, for any $x \in S$, $\partial \delta_{S} (x) =N_{S}(x)$. 

Some calculus rules for the limiting/Mordukhovich subdifferential used later are collected in the following lemmas (see \cite{Bot-Dao-Li,Mordukhovich2006,Mordukhovich2018,Rockafellar1998}).

\begin{lemma}[Fermat rule, see {\cite[Proposition 1.114]{Mordukhovich2006}}]
\label{lema22}
If a proper function $\varphi \colon \mathbb{R}^n \to (-\infty, +\infty]$ has a local minimum at $\bar{x},$ then $0 \in \widehat{\partial}\varphi(\bar{x})\subseteq \partial \varphi(\bar{x}).$	
\end{lemma}

\begin{lemma}[See {\cite[Theorem 4.10]{Mordukhovich2018}}]
\label{lema23}
Let $\varphi_i \colon \mathbb{R}^n \to (-\infty, +\infty]$, $i=1,\dots,m$ with $m \geq 2$, be locally Lipschitz around $\bar{x}$. Then the subdifferential of maximum function holds
\begin{align*}
\partial (\max \varphi_i) (\bar{x}) \subseteq \left\{    \sum_{i \in I(\bar{x})} \lambda_i \partial \varphi_i(\bar{x}) \mid \lambda_i \ge 0, \sum_{i \in I(\bar{x})} \lambda_i=1  \right\},
\end{align*}	
where $I(\bar{x}):=\{ i \in \{1,\ldots,m\}\mid \varphi_i(\bar{x})= (\max \varphi_i) (\bar{x})\}$.
\end{lemma}

We say that a function $\varphi\colon \mathbb{R}^n\to (-\infty, +\infty]$ is \emph{$\beta$-weakly convex on $S\subseteq \mathbb{R}^n$} if $\beta\geq 0$ and $\varphi +\delta_S +\frac{\beta}{2}\|\cdot\|^2$ is convex.

\begin{lemma}[See {\cite[Lemmas 2.2 and 4.1]{Bot-Dao-Li}}]
\label{l:weakcvx} 
Let $S$ be a nonempty subset of $\mathbb{R}^n$, let $\varphi\colon \mathbb{R}^n\to (-\infty, +\infty]$ be $\beta$-weakly convex on an open set containing $S$, and let $\bar x\in S$. Then the following hold:
\begin{enumerate}
\item 
$\varphi$ is continuous at $\bar x$ and $\partial \varphi(\bar x)$ is a nonempty closed convex subset of $\mathbb{R}^n$.
\item\label{l:weakcvx_inequality} 
For all $u\in\partial \varphi(\bar x)$ and all $x\in S$, $\varphi(x) -\varphi(\bar x) +\frac{\beta}{2}\|x -\bar x\|^2\geq \langle u, x -\bar x\rangle$.
\end{enumerate}
\end{lemma}

\section{Optimality conditions}
\label{s:optcond}

In this section, we derive and analyze necessary and sufficient conditions for the multiobjective optimization problem
\begin{align}\label{eq:prob}
\min_{x\in S} F(x),\tag{MP}
\end{align}
where $S$ is a nonempty closed subset of $\mathbb{R}^n$, $F\colon \mathbb{R}^n\to \mathbb{R}^m$ is a vector-valued function with $F :=(F_1, \dots, F_m)^\top$, and, for each $i\in I:=\{1, \dots, m\}$, 
\begin{align*}
F_i := f_i +g_i -h_i,    
\end{align*}
where $f_i\colon \mathbb{R}^n\to (-\infty, +\infty]$ is locally Lipschitz, $g_i\colon (-\infty, +\infty]$ is differentiable with an $\ell$-Lipschitz continuous gradient, and $h_i\colon \mathbb{R}^n\to \mathbb{R}$ is $\beta$-weakly convex, each satisfying its respective property on an open set containing $S$.

\begin{definition} 
Let $\bar x\in S$. We say that
\begin{enumerate}
\item 
$\bar x$ is a local weak Pareto solution of \eqref{eq:prob} if there exists a neighborhood $U$ of $\bar x$ such that there is no $x\in S\cap U$ satisfying  
\begin{align}\label{equa-1}
F_i(x) <F_i(\bar{x}),\quad i\in I.
\end{align}

\item 
$\bar x$ is a (global) weak Pareto solution of \eqref{eq:prob} if there is no $x\in S$ satisfying \eqref{equa-1}. 
\end{enumerate}
\end{definition}

\begin{remark}\label{Remark-1} 
By definition, it is easy to see that a point $\bar x\in S$ is a (local) weak Pareto solution of \eqref{eq:prob} iff $\bar x$ is a (local) solution of the following minimax programming problem
\begin{align*}
\min_{x\in S} \max_{i\in I} (F_i(x)-F_i(\bar x)). 
\end{align*}
\end{remark}

\begin{definition}
\label{d:stationary} 
Let $\bar x\in S$. We say that
\begin{enumerate}
\item 
$\bar x$ is a \emph{stationary point} of \eqref{eq:prob} if there exists $\lambda =(\lambda_1, \dots, \lambda_m)\in \mathbb{R}^m_+\setminus\{0\}$ such that
\begin{align*}
0\in \sum_{i\in I} \lambda_i\left(\partial f_i(\bar x) +\nabla g_i(\bar x) -\partial h_i(\bar x)\right) +N_S(\bar x).
\end{align*}

\item 
$\bar x$ is a \emph{strong stationary point} of \eqref{eq:prob} if there exists $\lambda =(\lambda_1, \dots, \lambda_m)\in \mathbb{R}^m_+\setminus\{0\}$ such that
\begin{align*}
\sum_{i\in I} \lambda_i\partial h_i(\bar x) \subseteq \sum_{i\in I} \lambda_i\left(\partial f_i(\bar x) +\nabla g_i(\bar x)\right) +N_S(\bar x).
\end{align*}
\end{enumerate}
\end{definition}

\begin{remark}
\begin{enumerate}
\item
Since $N_S(\bar{x})$ is a cone, the stationarity conditions in Definition~\ref{d:stationary} are positively homogeneous with respect to the multiplier vector. Hence, any $\lambda =(\lambda_1, \dots, \lambda_m)\in \mathbb{R}^m_+\setminus \{0\}$ satisfying either of stationarity conditions can be normalized so that $\sum_{i\in I} \lambda_i =1$. In the sequel, we may therefore assume without loss of generality that $\lambda\in \mathbb{R}^m_+$ with $\sum_{i\in I} \lambda_i =1$.
\item
Clearly, every strong stationary point of \eqref{eq:prob} is also a stationary one.  The reverse does not hold even for the scalar case; see, e.g., \cite[p.~133]{An-Nam-17}. However, it is easy to see that if the functions $h_i$, $i\in I$, are differentiable at $\bar x$ and $\bar x$ is a stationary point of \eqref{eq:prob}, then it is also a strong stationary point of \eqref{eq:prob}. 
\end{enumerate}
\end{remark}

\begin{theorem}[Necessary optimality condition]
\label{Necessary-Theorem}
Suppose that, for each $i\in I$, restricted to an open set containing $S$, $f_i\colon \mathbb{R}^n\to (-\infty, +\infty]$ is locally Lipschitz, $g_i\colon \mathbb{R}^n\to (-\infty, +\infty]$ is continuously differentiable, and $h_i\colon \mathbb{R}^n\to \mathbb{R}$ is $\beta$-weakly convex. If $\bar x$ is a local weak Pareto solution of \eqref{eq:prob}, then it is a stationary point of \eqref{eq:prob}.
\end{theorem}
\begin{proof} 
For each $i\in I$, let $u_i\in\partial h_i(\bar x)$ and let $\varphi_i\colon \mathbb{R}^n\to (-\infty, +\infty]$ be defined by
\begin{align*} 
\forall x\in \mathbb{R}^n,\quad \varphi_i(x) :=f_i(x) -f_i(\bar x) +g_i(x) -g_i(\bar x) -\langle u_i, x-\bar x\rangle +\frac{\beta}{2}\|x -\bar x\|^2.
\end{align*}
Then, for all $i\in I$, $\varphi_i(\bar x)=0$. Next, let $\varphi\colon \mathbb{R}^n\to (-\infty, +\infty]$ be defined by 
\begin{align*}
\forall x\in \mathbb{R}^n,\quad \varphi(x) :=\max_{i\in I} \varphi_i(x).    
\end{align*}
We claim that $\bar x$ is a local solution of the problem
\begin{align}\label{Aux-Problem}
\min_{x\in S} \varphi(x).
\end{align}
By Remark \ref{Remark-1} and the fact that $\bar x$ is a local weak Pareto solution of \eqref{eq:prob}, there exists a neighborhood $U$ of $\bar x$ such that
\begin{align*}
\forall x\in S\cap U,\quad \max_{i\in I} (F_i(x)-F_i(\bar x))\geq 0.
\end{align*}
This means that, for each $x\in S\cap U$, there exists $i_0\in I$ such that $F_{i_0}(x)-F_{i_0}(\bar x)\geq 0$,
or equivalently,
\begin{align*}
f_{i_0}(x) -f_{i_0}(\bar x) +g_{i_0}(x) -g_{i_0}(\bar x) -\left(h_{i_0}(x) -h_{i_0}(\bar x)\right) \geq 0.
\end{align*}
Since $u_{i_0}\in \partial h_{i_0}(\bar x)$, it follows from Lemma~\ref{l:weakcvx}\ref{l:weakcvx_inequality} that
\begin{align*}
h_{i_0}(x)-h_{i_0}(\bar x) +\frac{\beta}{2}\|x -\bar x\|^2\geq \langle u_{i_0}, x-\bar x\rangle. 
\end{align*}
Hence,
\begin{align*}
\varphi_{i_0}(x) &=f_{i_0}(x) -f_{i_0}(\bar x) +g_{i_0}(x) -g_{i_0}(\bar x) -\langle u_{i_0}, x-\bar x\rangle +\frac{\beta}{2}\|x -\bar x\|^2 \\
&\geq f_{i_0}(x) -f_{i_0}(\bar x) +g_{i_0}(x) -g_{i_0}(\bar x) -\left(h_{i_0}(x) -h_{i_0}(\bar x)\right) \geq 0.
\end{align*}
This implies that, for all $x\in S\cap U$, $\varphi(x)=\max_{i\in I}\varphi_i (x)\geq 0$. Thus, $\bar x$ is a local solution of \eqref{Aux-Problem}. By Lemma \ref{lema22}, one has 
\begin{align*}
0\in \partial(\varphi+\delta_S)(\bar x).
\end{align*} 
By assumptions, it is clear that $\varphi$ is locally Lipschitz and hence
\begin{align*}
0\in \partial\varphi(\bar x)+N_S(\bar x).
\end{align*} 
Thanks to Lemma \ref{lema23}, there exists $\lambda\in\mathbb{R}^m_+$ such that $\sum_{i\in I}\lambda_i=1$ and
\begin{align*}
0&\in \sum_{i\in I}\lambda_i\partial \varphi_i(\bar x)+N_S(\bar x)
\\
&=\sum_{i\in I}\lambda_i\left(\partial f_i(\bar x) +\nabla g_i(\bar x)\right) -\sum_{i\in I}\lambda_i u_i+N_S(\bar x)
\\
&\subseteq \sum_{i\in I}\lambda_i\left(\partial f_i(\bar x) +\nabla g_i(\bar x) -\partial h_i(\bar x)\right) +N_S(\bar x).
\end{align*}
This completes the proof of the theorem.
\end{proof}

\begin{remark} 
In \cite[Theorem 2.1]{Ji-Goh-Souza-16} (see also \cite[Theorem 3.1]{Qu-et al}), it is claimed that if, for each $i\in I$, $f_i$ and $h_i$ are convex, $g_i\equiv 0$, the constraint set $S$ is convex, and $\bar x$ is a local weak Pareto solution of \eqref{eq:prob}, then $\bar x$ is a strong stationary point of \eqref{eq:prob}. However, the following example shows that this assertion does not hold.
\end{remark}

\begin{example}\label{example-1} 
Let us consider problem \eqref{eq:prob} in the case where $S=\mathbb{R}$, $m=2$, and, for $i\in \{1, 2\}$, $g_i\equiv 0$, while $f_i$ and $h_i$ are convex functions defined by
\begin{align*}
\forall x\in\mathbb{R},\quad f_1(x)=|x|,\ f_2(x)=2x,\ h_1(x)=|x|+x,\ h_2(x)=|x|.
\end{align*}
Then  $F_1(x)=-x$ and $F_2(x)=2x-|x|$.  Let $\bar x=0$. We claim that $\bar x$ is a weak Pareto solution of \eqref{eq:prob}. Indeed,  one has
\begin{align*}
\forall x\in\mathbb{R},\quad \max\{F_1(x), F_2(x)\}=\max\{-x,2x-|x|\}\geq 0,
\end{align*}
as required. By Theorem \ref{Necessary-Theorem}, $\bar x$ is a stationary point of \eqref{eq:prob}.  We now show that $\bar x$ is not a strong stationary point of \eqref{eq:prob}. If otherwise, then there exist $\lambda_1\geq 0$, $\lambda_2\geq 0$ such that $\lambda_1+\lambda_2=1$ and
\begin{align}\label{equa-1-new}
\lambda_1\partial h_1(\bar x)+\lambda_2\partial h_2(\bar x)\subseteq \lambda_1 \partial f_1(\bar x)+\lambda_2\partial f_2(\bar x).
\end{align}
A direct computation shows that
\begin{align*}
&\partial f_1(\bar x)=[-1, 1],\ \partial f_2(\bar x)=\{2\},
\\
& \partial h_1(\bar x)=[0, 2], \text{~and~} \partial h_2(\bar x)=[-1, 1].
\end{align*}
Hence, \eqref{equa-1-new} means that
\begin{align*}
\lambda_1 [0,2]+\lambda_2 [-1, 1]\subseteq \lambda_1 [-1, 1]+\lambda_2.2,
\end{align*}
or, equivalently,
\begin{align*}
[-\lambda_2, 2\lambda_1+\lambda_2]\subseteq [-\lambda_1+2\lambda_2, \lambda_1+2\lambda_2].
\end{align*}
Thus,
\begin{align*}
\begin{cases}
-\lambda_1+2\lambda_2\leq  -\lambda_2,
\\
2\lambda_1+\lambda_2 \leq \lambda_1+2\lambda_2,
\\
\lambda_1\geq 0, \lambda_2\geq 0, \lambda_1+\lambda_2=1.
\end{cases}
\end{align*}
This implies that
\begin{align*}
\begin{cases} 
\lambda_2\leq \frac{1}{4},
\\
\lambda_2\geq \frac{1}{2},
\end{cases}
\end{align*}
a contradiction.
\end{example}

The following result gives a condition under which every local weak Pareto solution of \eqref{eq:prob} is a strong stationary point even when the functions $h_i$ are nonsmooth.

\begin{theorem} 
\label{t:strong-stationarity}
Suppose that the assumptions of Theorem \ref{Necessary-Theorem} are satisfied. If $\bar x$ is a local weak Pareto solution of \eqref{eq:prob} and there exist a neighborhood $U$ of $\bar x$ and a vector $  \lambda\in\mathbb{R}^m_+\setminus\{0\}$ such that
\begin{align}\label{eq:local-supportedness}
\forall x\in U\cap S,\quad \sum_{i \in I} \lambda _i F_{i}(x)\geq \sum_{i \in I} \lambda _i F_{i}(\bar x),
\end{align}
then $\bar x$ is a strong stationary point of \eqref{eq:prob}.
\end{theorem}	
\begin{proof} 
For each $i\in I$, let $u_i\in\partial h_i(\bar x)$ be arbitrary. Let $\varphi\colon \mathbb{R}^n\to (-\infty, +\infty]$ be the function defined by
\begin{align*}
\forall x\in \mathbb{R}^n,\quad \varphi(x) :=\sum_{i \in I} \lambda_i\left(f_i(x) +g_i(x) -h_i(\bar x) -\langle u_i, x -\bar x\rangle +\frac{\beta}{2}\|x -\bar x\|^2\right).
\end{align*}
In view of Lemma~\ref{l:weakcvx}\ref{l:weakcvx_inequality},
\begin{align*}
\forall x\in S,\quad h_i(x) -h_i(\bar x) +\frac{\beta}{2}\|x -\bar x\|^2\geq \langle u_i, x -\bar x\rangle,
\end{align*}
which together with \eqref{eq:local-supportedness} implies that
\begin{align*}
\forall x\in U\cap S,\quad \varphi(x) \geq \sum_{i \in I} \lambda_i (f_i(x) +g_i(x) -h_i(x)) =\sum_{i \in I} \lambda_i F_i(x) \geq \sum_{i \in I} \lambda_i F_i(\bar x) =\varphi(\bar x).
\end{align*}
This means that $\bar x$ is a local minimizer of $\varphi$ on $S$. Thus, by Lemma \ref{lema22} and the local Lipschitz continuity of $\varphi$, one has
\begin{align*}
0\in \partial\varphi(\bar x) +N_S(\bar x) =\sum_{i \in I} \lambda_i(\partial f_i(\bar x) +\nabla g_i(\bar x) -u_i) +N_S(\bar x),
\end{align*}  
or, equivalently,
\begin{align*}
\sum_{i \in I} \lambda_i u_i\in \sum_{i \in I} \lambda_i(\partial f_i(\bar x)+\nabla g_i(\bar x))+N_S(\bar x). 
\end{align*}
As each $u_i\in \partial h_i(\bar x)$ is chosen
arbitrarily, we obtain that $\bar x$ is a strong stationary point of \eqref{eq:prob}.
\end{proof}

\begin{remark}
By \cite[Theorem~3.5]{Ehrgott}, if the constraint set $S$ and the functions $F_i$, $i\in I$, are convex, then there exist a neighborhood $U$ of $\bar x$ and $\lambda\in\mathbb{R}^m_+\setminus\{0\}$ such that \eqref{eq:local-supportedness} holds. Hence, by Theorem~\ref{t:strong-stationarity}, every local weak Pareto solution in this convex setting is a strong stationary point.
\end{remark}

We now provide a sufficient condition for a strong stationary point of \eqref{eq:prob} to be a local weak Pareto solution of \eqref{eq:prob}.    

\begin{theorem}[Sufficient optimality condition, type~I]
\label{sufficient-theorem-typeI} 
Suppose that $S$ is convex and, for each $i\in I$, $f_i$ and $g_i$ are convex, and $h_i$ is polyhedral convex. If $\bar x\in S$ is a strong stationary point of \eqref{eq:prob}, then $\bar x$ is a local weak Pareto solution of \eqref{eq:prob}.
\end{theorem}
\begin{proof}
This follows from \cite[Example 3.1 and Corollary 3.4]{Gadhi-04}.    
\end{proof}

\begin{remark} 
It follows from Example \ref{example-1} that the converse conclusion  of Theorem \ref{sufficient-theorem-typeI} does not hold even for the case where the functions $f_i$ and $h_i$, $i\in I$, are polyhedral convex. The class of polyhedral convex functions includes finite maxima of affine functions and piecewise-affine convex functions, such as the $\ell_1$- and $\ell_\infty$-norms and hinge-type loss functions, including the standard hinge loss used in support vector machines.
\end{remark}

We next introduce a class of functions for which every stationary point of \eqref{eq:prob} is a local weak Pareto solution. To this end, we recall the concept of the contingent cone (also known as the Bouligand--Severi tangent cone) to the set $S$; see, e.g., \cite{Ehrgott,Jahn-04,Rockafellar1998,Sawaragi-et al-85}. The \emph{contingent cone} to $S$ at $\bar x\in S$ is defined by 
\begin{align*}
T_S(\bar x):=\{\nu\in\mathbb{R}^n\,|\, \exists t_k\to 0^+, \nu^k\to \nu\ \ \text{such that}\ \  \bar x+t_k\nu^k\in S\ \ \forall k\in\mathbb{N}\}.
\end{align*}
Note that when $S$ is convex, then $T_S(\bar x)=[N_S(\bar x)]^\circ$, where
\begin{align*}
[N_S(\bar x)]^\circ:=\{\nu\,|\, \langle \nu, u\rangle\leq 0 \ \ \forall u\in N_S(\bar x)\}.
\end{align*}

\begin{definition} 
We say that $(f, g, h)$ is \emph{locally generalized convex-concave} on $S$ at $\bar x\in S$ if there exists a neighborhood $U$ of $\bar x$ such that, for any $x\in U\cap S$, $u_i\in \partial f_i(\bar x)$, and $v_i\in\partial h_i(\bar x)$, $i\in I$, there exists $\nu\in T_S(\bar x)$ satisfying
\begin{align*}
&f_i(x)-f_i(\bar x)\geq \langle u_i, \nu\rangle,\quad  i\in I,
\\
&g_i(x)-g_i(\bar x)\geq \langle\nabla g_i(\bar x), \nu\rangle,\quad  i\in I,
\\
&h_i(\bar x)-h_i(x)\geq \langle -v_i, \nu\rangle,\quad  i\in I.
\end{align*}  
\end{definition}

\begin{remark}
\label{r:generalized convex-concave}
Suppose that $S$ is convex and let $\bar x\in S$. It is easy to see that if $f_i$ and $g_i$, $i\in I$, are convex and there exists a neighborhood $U$ of $\bar x$ such that
\begin{align*}
\forall x\in U,\quad h_i(x)=\langle a_i, x\rangle+b_i,\quad i\in I,    
\end{align*}
then $(f, g, h)$ is locally generalized convex-concave on $S$ at $\bar x$ with $\nu =x -\bar x$. In particular, this situation arises locally when $h_i$ is piecewise affine and $\bar x$ lies in a region on which a single affine piece is active.
\end{remark}

The following result gives a sufficient condition for a stationary point of \eqref{eq:prob} to be a local weak Pareto solution.

\begin{theorem}[Sufficient optimality condition, type~II]
\label{sufficient-theorem} 
Suppose that $S$ is convex and $(f, g, h)$ is locally generalized convex-concave on $S$ at $\bar x\in S$. If $\bar x$ is a stationary point of  \eqref{eq:prob}, then it is a local weak Pareto solution of \eqref{eq:prob}.
\end{theorem}
\begin{proof} Since $\bar x$ is a stationary point of \eqref{eq:prob}, there exist $\lambda_i\geq 0$, $i\in I$, such that $\sum_{i\in I}\lambda_i=1$ and
\begin{align*}
0\in \sum_{i\in I} \lambda_i(\partial f_i(\bar x)+\nabla g_i(\bar x)-\partial h_i(\bar x)) +N_S(\bar x).
\end{align*} 
Hence, there are $u_i\in\partial f_i(\bar x)$, $v_i\in\partial h_i(\bar x)$, and $w\in N_S(\bar x)$ such that
\begin{align}\label{s-stationary-1}
\sum_{i\in I}\lambda_i(u_i+\nabla g_i(\bar x)-v_i)=-w.
\end{align}
Since $(f, g, h)$ is locally generalized convex-concave on $S$ at $\bar x$, there  exists a neighborhood $U$ of $\bar x$ such that for any $x\in U\cap S$, we find $\nu\in T_S(\bar x)$  satisfying
\begin{align*}
&f_i(x)-f_i(\bar x)\geq \langle u_i, \nu\rangle,\quad  i\in I,
\\
&g_i(x)-g_i(\bar x)\geq \langle\nabla g_i(\bar x), \nu\rangle,\quad  i\in I,
\\
&h_i(\bar x)-h_i(x)\geq \langle -v_i, \nu\rangle,\quad  i\in I.
\end{align*}   
Combining the above inequalities together with \eqref{s-stationary-1} we arrive at
\begin{align*}
&\sum_{i\in I} \lambda_i(f_i(x)+g_i(x)-h_i(x)) -\sum_{i\in I} \lambda_i(f_i(\bar x)+g_i(\bar x)-h_i(\bar x)) \\
&\geq \sum_{i\in I} \lambda_i(\langle u_i, \nu\rangle +\langle\nabla g_i(\bar x), \nu\rangle -\langle v_i, \nu\rangle)
\\
&=\left\langle -w, \nu\right\rangle\geq 0,
\end{align*}
which implies that
\begin{align*}
\forall x\in U\cap S,\quad \max_{i\in I} (F_i(x)-F_i(\bar x))\geq 0.
\end{align*}
Thus, $\bar x$ is a local weak Pareto solution of \eqref{eq:prob}.
\end{proof}

\section{Proposed Algorithm}
\label{s:algo}

In this section, we propose and analyze an algorithm for solving \eqref{eq:prob} that exploits the composite DC-type structure of the objective functions. At each iteration, the update is obtained by minimizing over the constraint set $S$ the maximum of local models of the componentwise objective variations together with a proximal regularization term. More precisely, the possibly nonsmooth term $f_i$ is retained in the model, while the smooth term $g_i$ and the weakly convex subtractive term $h_i$ are linearized at the current iterate using the gradient and a subgradient, respectively. No weight vector for scalarizing the original vector objective is fixed in advance, and possible weights arise only implicitly through the active component models in the optimality condition of the max-type subproblem. Related max-type constructions have been used in proximal-gradient methods for multiobjective optimization; see, e.g., \cite{Tanabe-Fukuda-Yamashita-19}. Algorithm~\ref{algo:PSG} extends this construction to constrained DC-type settings, where the constraint set $S$ is incorporated directly, $f_i$ is locally Lipschitz and possibly nonconvex, and $h_i$ is weakly convex. This extension admits a simultaneous descent estimate for all objective functions, which plays a central role in the convergence analysis.

\begin{tcolorbox}[
left=0pt,right=0pt,top=0pt,bottom=0pt,
colback=blue!10!white, colframe=blue!50!white,
boxrule=0.2pt,
breakable]
\begin{algorithm}
\label{algo:PSG}
\step{}
Let $x^0 \in S$ and $\bar{\gamma}\in \left(0, \frac{1}{\ell +\beta}\right)$, where we adopt the convention that $\frac{1}{\ell +\beta} =+\infty$ if $\ell +\beta =0$. Set $k =0$.

\step{}\label{step:PSG-main}
For each $i\in I$, let $u_i^k\in \partial h_i(x^k)$. Let $\gamma_k\in (0, \bar{\gamma}]$ and find
\begin{align*}
x^{k+1} \in \argmin_{x\in S} \left(\max_{i\in I} \left(f_i(x) -f_i(x^k) +\langle \nabla g_i(x^k) -u_i^k, x -x^k\rangle\right) +\frac{1}{2\gamma_k}\|x -x^k\|^2\right).
\end{align*}

\step{}
If a termination criterion is not met, set $k =k +1$ and go to Step~\ref{step:PSG-main}.
\end{algorithm}
\end{tcolorbox}

\begin{remark}
In the case when all $h_i\equiv 0$, Step~\ref{step:PSG-main} of Algorithm~\ref{algo:PSG} reduces to
\begin{align*}
x^{k+1} \in \argmin_{x\in S} \left(\max_{i\in I} \left(f_i(x) -f_i(x^k) +\langle \nabla g_i(x^k), x -x^k\rangle\right) +\frac{1}{2\gamma_k}\|x -x^k\|^2\right),
\end{align*}
which, up to notation and parameterization, is the subproblem used in \cite[Algorithm~3.2]{Tanabe-Fukuda-Yamashita-19} if $S =\mathbb{R}^n$ and every $f_i$ is convex. If, instead, all $g_i\equiv 0$, then Step~\ref{step:PSG-main} becomes
\begin{align*}
x^{k+1} \in \argmin_{x\in S} \left(\max_{i\in I} \left(f_i(x) -f_i(x^k) -\langle u_i^k, x -x^k\rangle\right) +\frac{1}{2\gamma_k}\|x -x^k\|^2\right),
\end{align*}
which is closely related to the subproblem in \cite[Algorithm~1]{Ji-Goh-Souza-16}. However, the latter assumes that $f_i$ and $h_i$ are convex and uses the uncentered models $f_i(x) -\langle u_i^k, x -x^k\rangle$. Thus, the two subproblems are generally different because the centering terms $f_i(x^k)$ depend on the objective index $i$. In general, Step~\ref{step:PSG-main} requires solving a subproblem of a max structure that admits several equivalent reformulations, which can be exploited for numerical solution.

\begin{enumerate}
\item
\emph{(Epigraph reformulation)} Introducing an auxiliary scalar variable $t\in\mathbb{R}$, the subproblem can be formulated as
\begin{align*}
\min_{x\in S,\; t\in\mathbb{R}}~~ &t +\frac{1}{2\gamma_k}\|x -x^k\|^2 \\
\text{subject to}~~ &f_i(x) -f_i(x^k) +\langle \nabla g_i(x^k) -u_i^k, x -x^k\rangle \leq t,\quad  i\in I.
\end{align*}
This formulation separates the pointwise maximum from the quadratic regularization term and reduces Step~\ref{step:PSG-main} to a single-objective constrained optimization problem. From an implementation viewpoint, the latter can be handled using standard nonlinear optimization techniques.

\item
\emph{(Simplex-weight representation)} Using the identity
\begin{align*}
\max_{i\in I} a_i =\max_{\lambda\in\Delta} \sum_{i\in I}\lambda_i a_i,
\quad
\Delta :=\Big\{\lambda\in \mathbb{R}_+^m \mid \sum_{i\in I} \lambda_i =1\Big\},
\end{align*}
the subproblem admits the equivalent formulation
\begin{align*}
\min_{x\in S} \max_{\lambda\in\Delta} \left(\sum_{i\in I} \lambda_i\left(f_i(x) -f_i(x^k) +\langle \nabla g_i(x^k) -u_i^k, x -x^k\rangle\right) +\frac{1}{2\gamma_k}\|x -x^k\|^2\right).
\end{align*}
This representation shows that Step~\ref{step:PSG-main} can be interpreted as computing a proximal-type step associated with a weighted aggregation of the objective variations, where the weights $\lambda_i$ belong to the simplex and are determined implicitly by the subproblem. Although this does not necessarily simplify the numerical solution of Step~\ref{step:PSG-main}, it is useful for interpreting the multipliers associated with the active component models and can be convenient when exploiting the structure of particular problems.
\end{enumerate}
\end{remark}

We now establish the main convergence properties of Algorithm~\ref{algo:PSG}, including descent estimates, boundedness of the generated sequence, and stationarity of its cluster points.

\begin{theorem}
\label{t:cvg}
Let $(x^k)_{k\in \mathbb{N}}$ be a sequence generated by Algorithm~\ref{algo:PSG}. Then the following hold:
\begin{enumerate}
\item\label{t:cvg_decrease}
For all $i\in I$ and all $k\in \mathbb{N}$,
\begin{align*}
F_i(x^{k+1}) +\left(\frac{1}{2\gamma_k} -\frac{\ell +\beta}{2}\right)\|x^{k+1} -x^k\|^2\leq F_i(x^k).    
\end{align*}

\item\label{t:cvg_asymptotic}
If there exists $i\in I$ such that $F_i$ is bounded from below on $S$, then
\begin{align*}
\sum_{k=0}^{+\infty} \|x^{k+1} -x^k\|^2 <+\infty    
\end{align*}
and consequently, $x^{k+1} -x^k\to 0$ as $k\to +\infty$.

\item\label{t:cvg_bounded}
If there exists $i\in I$ such that the set $\{x\in S: F_i(x)\leq F_i(x^0)\}$ is bounded, then the sequence $(x^k)_{k\in \mathbb{N}}$ is bounded.

\item\label{t:cvg_crit}
If $\liminf_{k\to +\infty} \gamma_k >0$ and there exists $i\in I$ such that $F_i$ is bounded from below on $S$, then every cluster point of the sequence $(x^k)_{k\in \mathbb{N}}$ is a stationary point of \eqref{eq:prob}.
\end{enumerate}
\end{theorem}
\begin{proof}
\ref{t:cvg_decrease}: We first see that, for all $k\in \mathbb{N}$, $x^k\in S$. Let $i\in I$ and $k\in \mathbb{N}$. It follows from Step~\ref{step:PSG-main} of Algorithm~\ref{algo:PSG} that, for all $x\in S$,
\begin{align*}
&f_i(x^{k+1}) -f_i(x^k) +\langle \nabla g_i(x^k) -u_i^k, x^{k+1} -x^k\rangle +\frac{1}{2\gamma_k}\|x^{k+1} -x^k\|^2 \\
&\leq \max_{i\in I} \left(f_i(x^{k+1}) -f_i(x^k) +\langle \nabla g_i(x^k) -u_i^k, x^{k+1} -x^k\rangle\right) +\frac{1}{2\gamma_k}\|x^{k+1} -x^k\|^2 \\
&\leq \max_{i\in I} \left(f_i(x) -f_i(x^k) +\langle \nabla g_i(x^k) -u_i^k, x -x^k\rangle\right) +\frac{1}{2\gamma_k}\|x -x^k\|^2.
\end{align*}
Letting $x =x^k$ yields 
\begin{align}\label{eq:fifi}
f_i(x^{k+1}) -f_i(x^k) +\langle \nabla g_i(x^k) -u_i^k, x^{k+1} -x^k\rangle +\frac{1}{2\gamma_k}\|x^{k+1} -x^k\|^2 \leq 0.
\end{align}

On the other hand, we have from
the Lipschitz continuity of $\nabla g_i$ and \cite[Lemma~1.2.3]{Nes18} that 
\begin{align*}
g_i(x^{k+1}) -g_i(x^k) \leq \langle \nabla g_i(x^k), x^{k+1} -x^k\rangle +\frac{\ell}{2}\|x^{k+1} -x^k\|^2,    
\end{align*}
and from the weak convexity of $h_i$ that
\begin{align*}
-h_i(x^{k+1}) +h_i(x^k) \leq -\langle u_i^k, x^{k+1} -x^k\rangle +\frac{\beta}{2}\|x^{k+1} -x^k\|^2.  
\end{align*}
Combining with \eqref{eq:fifi}, we obtain that
\begin{align*}
f_i(x^{k+1}) -f_i(x^k) +\frac{1}{2\gamma_k}\|x^{k+1} -x^k\|^2 +g_i(x^{k+1}) -g_i(x^k) -h_i(x^{k+1}) +h_i(x^k)\leq \frac{\ell +\beta}{2}\|x^{k+1} -x^k\|^2,   
\end{align*}
or equivalently,
\begin{align*}
F_i(x^{k+1}) +\left(\frac{1}{2\gamma_k} -\frac{\ell +\beta}{2}\right)\|x^{k+1} -x^k\|^2\leq F_i(x^k).    
\end{align*}

\ref{t:cvg_asymptotic}: Assume that there exists $i\in I$ such that $F_i$ is bounded from below on $S$. By combining with \ref{t:cvg_decrease}, the sequence $(F_i(x^k))_{k\in \mathbb{N}}$ is nonincreasing and bounded from below, so it is convergent. Now, telescoping the inequality in \ref{t:cvg_decrease} and noting that, for all $k\in \mathbb{N}$, $\gamma_k\in (0, \bar{\gamma}]$, we obtain
\begin{align*}
\left(\frac{1}{2\bar{\gamma}} -\frac{\ell +\beta}{2}\right)\sum_{k=0}^{+\infty} \|x^{k+1} -x^k\|^2\leq F_i(x^0) -\lim_{n\to +\infty} F_i(x^{n+1}) <+\infty.     
\end{align*}
Since $\bar{\gamma}\in \left(0, \frac{1}{\ell +\beta}\right)$, it follows that $\sum_{k=0}^{+\infty} \|x^{k+1} -x^k\|^2 <+\infty$, and hence $x^{k+1} -x^k\to 0$ as $k\to +\infty$.

\ref{t:cvg_bounded}: Assume that there exists $i\in I$ such that the set $S_0 :=\{x\in S: F_i(x)\leq F_i(x^0)\}$ is bounded. We derive from \ref{t:cvg_decrease} that, for all $k\in \mathbb{N}$, $F_i(x^k)\leq F_i(x^0)$, and so $x^k\in S_0$. The conclusion then follows.

\ref{t:cvg_crit}: Let $k\in \mathbb{N}$. We have from the optimality condition for the $x$-update that
\begin{align}\label{eq:optcond}
0 &\in \partial \left(\max_{i\in I} \left(f_i(\cdot) -f_i(x^k) +\langle \nabla g_i(x^k) -u_i^k, \cdot -x^k\rangle\right) +\frac{1}{2\gamma_k}\|\cdot -x^k\|^2+\delta_S(\cdot)\right)(x^{k+1}) \notag \\
&\subseteq \sum_{i\in I} \lambda_i^k\left(\partial f_i(x^{k+1}) +\nabla g_i(x^k) -u_i^k\right) +\frac{1}{\gamma_k}(x^{k+1} -x^k)+N_S(x^{k+1}),
\end{align}
where $\lambda_1^k, \dots, \lambda_m^k\in [0, +\infty)$ summing up to $1$.

Let $\bar{x}$ be a cluster point of the sequence $(x^k)_{k\in \mathbb{N}}$. Then there exists a subsequence $(x^{n_k})_{k\in \mathbb{N}}$ such that $x^{n_k}\to \bar{x}$ as $k\to +\infty$. Since $h_i +\frac{\beta}{2}\|\cdot\|^2$ is finite and convex on an open neighborhood of $S$, $\partial h_i =\partial(h_i +\frac{\beta}{2}\|\cdot\|^2) -\beta\Id$ is locally bounded around $\bar x$, and so  $(u_i^k)_{k\in \mathbb{N}}$ is bounded. By passing to a subsequence if necessary, we can assume that, for each $i\in I$, $u_i^{n_k}\to \bar{u}_i$  and $\lambda_i^{n_k}\to \lambda_i\in [0, +\infty)$ as $k\to +\infty$. Clearly, $\sum_{i=1}^m \lambda_i =1$ and by the robustness of the convex subdifferential, one has  $\bar u_i \in \partial h_i(\bar{x})$.  By \ref{t:cvg_asymptotic}, $x^{n_k+1}\to \bar{x}$. Replacing $k$ in \eqref{eq:optcond} by $n_k$, letting $k\to +\infty$, and using the robustness property of $\partial f_i$ and the continuity of $\nabla g_i$, we derive that
\begin{align*}
0 \in \sum_{i\in I} \lambda_i\left(\partial f_i (\bar{x}) +\nabla g_i(\bar x) -\bar{u}_i\right) + N_S(\bar x) \subseteq \sum_{i\in I} \lambda_i\left(\partial f_i (\bar{x}) +\nabla g_i(\bar x) -\partial h_i (\bar{x})\right) + N_S(\bar x) ,   
\end{align*}
which completes the proof.
\end{proof}

\section{Numerical experiments}
\label{s:numerical}

In this section, we illustrate the performance of Algorithm~\ref{algo:PSG} on a multi-sensor sparse recovery problem. An unknown sparse signal $x^\dagger\in\mathbb{R}^n$ is to be recovered from observations obtained by two measurement systems according to
\begin{align*}
b_i =A_i x^\dagger +\varepsilon_i,\quad i\in \{1, 2\},
\end{align*}
where $A_i\in \mathbb{R}^{p_i\times n}$ is the sensing matrix, $b_i\in \mathbb{R}^{p_i}$ is the observation vector, and $\varepsilon_i\in \mathbb{R}^{p_i}$ is the measurement noise. Since the two systems may have different noise levels and sensing characteristics, we treat the corresponding data-fidelity terms as separate objectives rather than combining them using prescribed scalarization weights. We consider
\begin{align}\label{eq:sparse-recovery}
\min_{x\in S}\ (F_1(x), F_2(x))^\top,\quad
F_i(x) :=\frac{\alpha_i}{2}\|A_i x -b_i\|^2 +\mu(\|x\|_1 -\|x\|),\quad 
i \in \{1, 2\},
\end{align}
where $S =[-M, M]^n$ and $\alpha_i :=\frac{1}{\|b_i\|^2}$. The normalization by $\|b_i\|^2$ makes the two data-fidelity terms comparable in scale. The parameter $\mu >0$ controls the weight of the nonconvex regularizer $\|\cdot\|_1 -\|\cdot\|$, which promotes sparsity and reduces the shrinkage effect associated with the convex $\|\cdot\|_1$ penalty. Problem~\eqref{eq:sparse-recovery} fits the framework of \eqref{eq:prob} with
\begin{align*}
f_i(x) =\mu\|x\|_1,\quad
g_i(x) =\frac{\alpha_i}{2}\|A_i x -b_i\|^2,\quad
h_i(x) =\mu\|x\|,\quad 
i \in \{1, 2\}.
\end{align*}
In particular, $\nabla g_i$ is $\ell$-Lipschitz continuous with $\ell :=\max\{\alpha_1\|A_1\|_2^2, \alpha_2\|A_2\|_2^2\}$, whereas $h_i$ is convex, and hence $\beta =0$.

We compare Algorithm~\ref{algo:PSG} with the proximal point algorithm proposed by Ji, Goh, and de~Souza (JGS-PPA) in \cite[Algorithm~1]{Ji-Goh-Souza-16}. To apply JGS-PPA, we write
\begin{align*}
F_i(x) =\widetilde{f}_i(x) -h_i(x),\quad
\widetilde{f}_i(x) :=\mu\|x\|_1 +g_i(x),\quad
h_i(x) :=\mu\|x\|,\quad 
i\in \{1,2\}.
\end{align*}
Since both algorithms require a subgradient of $h_i$ at the current iterate and $h_1 =h_2$, we use the common choice $u^k\in \partial h_i(x^k)$ with
\begin{align*}
u^k =
\begin{cases}
\mu\frac{x^k}{\|x^k\|} &\text{if~} x^k\neq 0, \\
0 &\text{if~} x^k =0,
\end{cases}
\quad i\in \{1,2\},
\end{align*}
and set
\begin{align*}
d_i^k :=\nabla g_i(x^k) -u^k,\quad i\in \{1,2\}.
\end{align*}
Since $f_1 =f_2 =\mu\|\cdot\|_1$, Step~\ref{step:PSG-main} of Algorithm~\ref{algo:PSG} becomes
\begin{align}\label{eq:PSG-sparse-subproblem}
x^{k+1}\in \argmin_{x\in S}\max_{\theta\in [0,1]} \left(\mu\|x\|_1 +\langle \theta d_1^k +(1-\theta)d_2^k, x -x^k\rangle +\frac{1}{2\gamma}\|x -x^k\|^2\right),
\end{align}
where we have used $\max\{a_1, a_2\} =\max_{\theta\in [0,1]} \{\theta a_1 +(1-\theta)a_2\}$. On the other hand, the subproblem of JGS-PPA can be written as
\begin{align}\label{eq:JGS-subproblem}
x^{k+1}\in \argmin_{x\in S}\max_{\eta\in [0,1]} \left(\mu\|x\|_1 +\eta g_1(x) +(1-\eta)g_2(x) -\langle u^k, x -x^k\rangle +\frac{1}{2\gamma}\|x -x^k\|^2\right).
\end{align}
Thus, the main computational difference is that Algorithm~\ref{algo:PSG} linearizes the smooth terms $g_i$, whereas JGS-PPA retains them in the subproblem.

We next describe how \eqref{eq:PSG-sparse-subproblem} is solved. Since $S =[-M, M]^n$ is compact and convex, and the objective in \eqref{eq:PSG-sparse-subproblem} is convex in $x$ and affine in $\theta$, the corresponding min-max problem admits the usual minimax interchange. For $\theta\in [0,1]$, consider the value function
\begin{align*}
\phi_k(\theta) :=\min_{x\in S} \left(\mu\|x\|_1 +\langle \theta d_1^k +(1-\theta)d_2^k, x -x^k\rangle +\frac{1}{2\gamma}\|x -x^k\|^2\right).
\end{align*}
Hence, solving \eqref{eq:PSG-sparse-subproblem}
amounts to maximizing the concave function $\phi_k$ over $[0,1]$.
For every fixed $\theta$, the minimizer in the definition of
$\phi_k$ is unique and is given by
\begin{align}\label{eq:x-theta}
x^k(\theta) =P_S\!\left(\prox_{\gamma\mu\|\cdot\|_1}(x^k-\gamma(\theta d_1^k +(1-\theta)d_2^k))\right),
\end{align}
where $\prox_{\tau\|\cdot\|_1}$ is the \emph{shrinkage}
operator given componentwise by
\begin{align*}
\left(\prox_{\tau\|\cdot\|_1}(z)\right)_j :=\operatorname{sign}(z_j)\max\{|z_j| -\tau, 0\}
\end{align*}
(see, e.g., \cite{Beck-Teboulle-09}) and $P_S$ denotes the projection onto $S$, i.e., coordinatewise clipping to $[-M, M]$.

Since the minimizer $x^k(\theta)$ is unique, Danskin's theorem gives
\begin{align*}
\phi_k'(\theta) =\langle d_1^k -d_2^k, x^k(\theta) -x^k\rangle.
\end{align*}
Since $\phi_k$ is concave, $\phi_k'$ is nonincreasing on $[0,1]$.
Thus, a maximizer $\theta_k$ of $\phi_k$ is characterized by
\begin{align}\label{eq:theta-k}
\theta_k=
\begin{cases}
0 &\text{if~} \phi_k'(0)\leq 0, \\
1 &\text{if~} \phi_k'(1)\geq 0, \\
\text{a root of~} \phi_k'(\theta) =0
&\text{if~} \phi_k'(0) >0 >\phi_k'(1),
\end{cases}
\end{align}
and $x^{k+1}=x^k(\theta_k)$. In the last case, the scalar equation is solved using Brent's method~\cite{Brent-73}. Consequently, each iteration of Algorithm~\ref{algo:PSG} requires only gradient evaluations, shrinkage operator, box projection, and, when necessary, a one-dimensional root search.

A similar scalarization applies to \eqref{eq:JGS-subproblem}. For $\eta\in [0,1]$, define
\begin{align}\label{eq:JGS-weighted}
\psi_k(\eta) :=\min_{x\in S} \left(\mu\|x\|_1 +\eta g_1(x) +(1-\eta)g_2(x) -\langle u^k, x -x^k\rangle +\frac{1}{2\gamma}\|x -x^k\|^2\right)
\end{align}
and denote its unique minimizer by $x_J^k(\eta)$. By the same minimax argument as above, solving \eqref{eq:JGS-subproblem} amounts to maximizing the concave function $\psi_k$ over $[0,1]$. Unlike \eqref{eq:x-theta}, $x_J^k(\eta)$ does not in general have a closed form because the quadratic data-fidelity terms $g_i$ are retained rather than linearized. We compute it by FISTA~\cite{Beck-Teboulle-09}. The smooth part of \eqref{eq:JGS-weighted} has Lipschitz constant
\begin{align*}
L_\eta :=\eta\alpha_1\|A_1\|_2^2 +(1-\eta)\alpha_2\|A_2\|_2^2 +\frac{1}{\gamma}
\end{align*}
and each FISTA proximal step is obtained by shrinkage operator followed by projection onto $S$.

Since $x_J^k(\eta)$ is unique, Danskin's theorem yields
\begin{align*}
\psi_k'(\eta) =g_1(x_J^k(\eta)) -g_2(x_J^k(\eta)).
\end{align*}
Since $\psi_k$ is concave, $\psi_k'$ is nonincreasing, and a maximizing $\eta_k$ is obtained from the same endpoint/root rule as in \eqref{eq:theta-k}, with $\phi_k'$ replaced by $\psi_k'$. Thus, JGS-PPA also involves a one-dimensional search, but each evaluation requires solving the inner strongly convex problem \eqref{eq:JGS-weighted} by FISTA, whereas Algorithm~\ref{algo:PSG} computes the corresponding candidate explicitly from \eqref{eq:x-theta}.

\paragraph{Experimental setup.} We generate $20$ independent test instances using random seeds $100, \dots, 119$. The dimensions and parameters are
\begin{align*}
n =128,\quad
p_1 =p_2 =40,\quad
\|x^\dagger\|_0 =10,\quad
M =3,\quad
\mu =3\times10^{-3}.
\end{align*}
The support of $x^\dagger$ is selected uniformly at random, while its nonzero entries have independent random signs and magnitudes uniformly distributed on $[1, 2]$. To generate sensing matrices of different difficulty, let $Z_i\in\mathbb{R}^{p_i\times n}$ have independent standard Gaussian entries and let $C_i\in\mathbb{R}^{n\times n}$ be the Toeplitz correlation matrix
\begin{align*}
(C_i)_{rs} =\rho_i^{|r-s|},\quad
\rho_1 =0.70,\quad
\rho_2 =0.90.
\end{align*}
We generate $A_i$ from $Z_iC_i^{1/2}$ and then normalize its columns to have unit Euclidean norm. The noise vectors are independently generated and scaled so that
\begin{align*}
\frac{\|\varepsilon_1\|}{\|A_1x^\dagger\|} =0.02,\quad
\frac{\|\varepsilon_2\|}{\|A_2x^\dagger\|} =0.05.
\end{align*}

To provide a common stable initialization, we initialize both algorithms from the solution of the convex $\ell_1$-regularized minimax problem obtained from \eqref{eq:sparse-recovery} by omitting the term $-\mu\|x\|$. The cost of computing this warm start is excluded from the reported CPU times. We use the constant stepsize $\gamma =0.95/\ell$. The iterations are terminated when
\begin{align*}
\frac{\|x^{k +1} -x^k\|}{\max\{1, \|x^k\|\}}\leq 10^{-5}.
\end{align*}
For JGS-PPA, each inner FISTA iteration is terminated when the corresponding relative step is at most $10^{-6}$, and the one-dimensional scalar searches are performed with tolerance $10^{-6}$.

\paragraph{Comparison of the results.} Table~\ref{tab:sparse-results} reports the mean and standard deviation over the $20$ instances. We use the relative reconstruction error
\begin{align*}
E(x) :=\frac{\|x -x^\dagger\|}{\|x^\dagger\|},
\end{align*}
the worst normalized sensor residual
\begin{align*}
R_{\max}(x) :=\max_{i\in \{1, 2\}} \frac{\|A_i x -b_i\|}{\|b_i\|},
\end{align*}
and the worst objective value
\begin{align*}
F_{\max}(x) :=\max\{F_1(x), F_2(x)\}.
\end{align*}
To assess support recovery, let
\begin{align*}
\mathcal{S}^\dagger :=\{j: x_j^\dagger\neq 0\},
\quad
\widehat{\mathcal{S}}(x) :=\left\{j: |x_j|\geq 0.05\max_{\ell}|x_\ell|\right\}.
\end{align*}
The support precision and recall are defined by
\begin{align*}
\operatorname{Prec}(x) :=\frac{|\widehat{\mathcal{S}}(x)\cap\mathcal{S}^\dagger|}{|\widehat{\mathcal{S}}(x)|},\quad
\operatorname{Rec}(x) :=\frac{|\widehat{\mathcal{S}}(x)\cap\mathcal{S}^\dagger|}{|\mathcal{S}^\dagger|},
\end{align*}
and the support recovery F-score is
\begin{align*}
\operatorname{Fscore}(x) :=\frac{2\operatorname{Prec}(x)\operatorname{Rec}(x)}{\operatorname{Prec}(x) +\operatorname{Rec}(x)}.
\end{align*}
We record exact support recovery when $\widehat{\mathcal{S}}(x) =\mathcal{S}^\dagger$. Next, observe that a nonzero $x\in S$ is a stationary point of \eqref{eq:sparse-recovery} if and only if there exists $\theta\in [0, 1]$ such that
\begin{align*}
0\in \partial(\mu\|\cdot\|_1 +\delta_S)(x) +\theta\nabla g_1(x) +(1-\theta)\nabla g_2(x) -\mu\frac{x}{\|x\|}.
\end{align*}
Using the proximal characterization of this inclusion with stepsize $1/\ell$, together with the fact that $S =[-M, M]^n$, we define the algorithm-independent stationarity residual
\begin{align*}
\mathcal{R}(x) :=\min_{\theta\in [0, 1]} \ell\left\|x -P_S\left(\prox_{\frac{\mu}{\ell}\|\cdot\|_1}\left(x -\frac{1}{\ell}\left[\theta\nabla g_1(x) +(1-\theta)\nabla g_2(x) -\mu\frac{x}{\|x\|}\right]\right)\right)\right\|,
\end{align*}
which vanishes precisely at stationary points.

\begin{table}[!htbp]
\centering
\caption{Comparison of Algorithm~\ref{algo:PSG} and JGS-PPA over $20$ test instances. Values are mean $\pm$ standard deviation unless otherwise indicated.}
\label{tab:sparse-results}
\begin{tabular}{lcc}
\hline
Metric & Algorithm~\ref{algo:PSG} & JGS-PPA\\
\hline
Relative reconstruction error $E$
& $0.08252\pm0.02009$
& $0.08252\pm0.02008$\\
Worst sensor residual $R_{\max}$
& $0.05653\pm0.00338$
& $0.05652\pm0.00339$\\
Worst objective value $F_{\max}$
& $0.030810\pm0.001594$
& $0.030810\pm0.001594$\\
support recovery F-score
& $0.9929\pm0.0174$
& $0.9929\pm0.0174$\\
Exact support recovery
& $17/20$
& $17/20$\\
Stationarity residual $\mathcal{R}$
& $(5.13\pm1.85)\times10^{-5}$
& $(5.16\pm1.86)\times10^{-5}$\\
Outer iterations
& $683.3\pm312.6$
& $685.2\pm313.1$\\
CPU time (s)
& $\mathbf{0.196\pm0.092}$
& $2.244\pm0.942$\\
\hline
\end{tabular}
\end{table}

\begin{figure}[!htbp]
\centering
\includegraphics[width=0.6\textwidth]{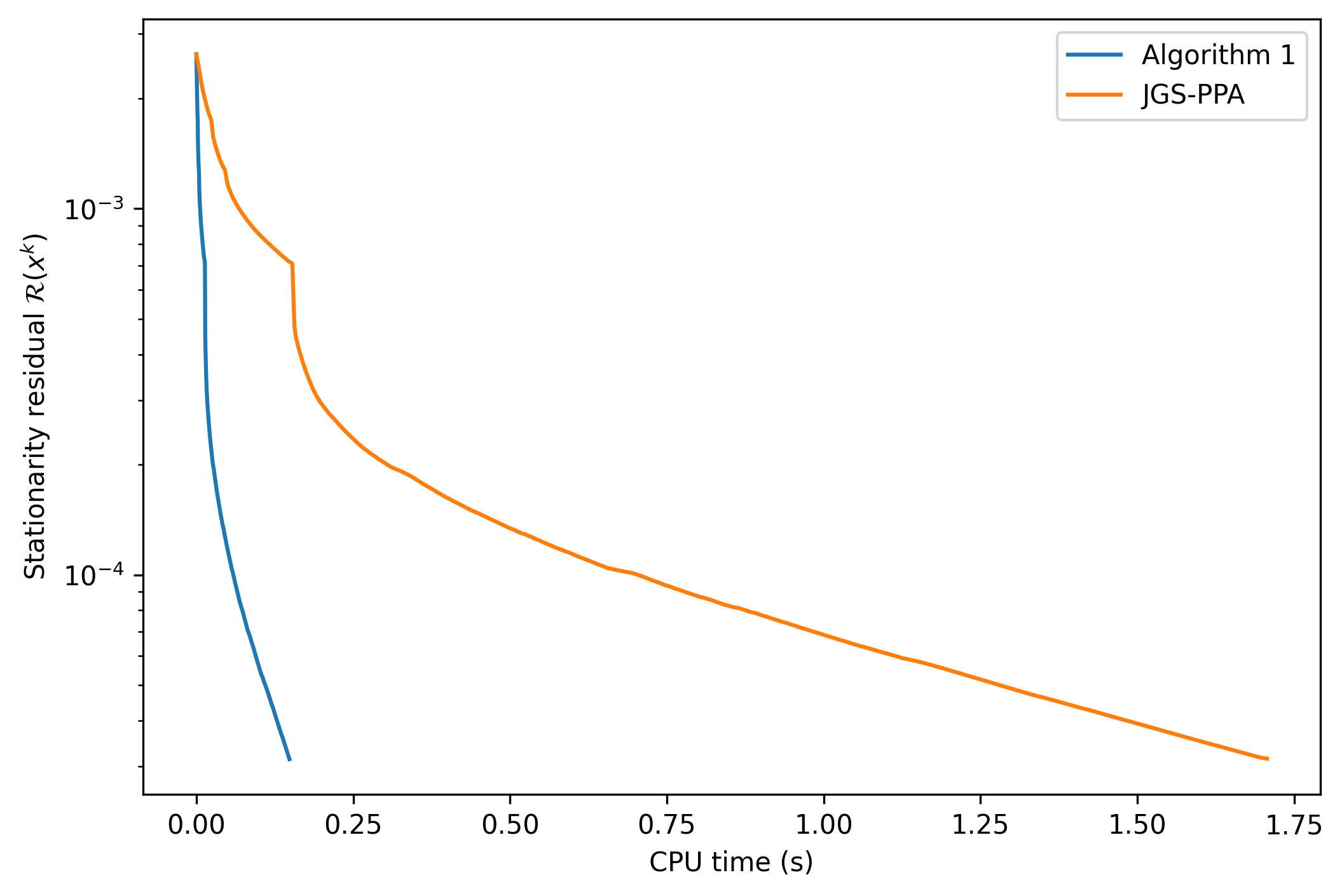}
\caption{Stationarity residual versus CPU time for a representative test instance.}
\label{fig:sparse-residual-time}
\end{figure}

As shown in Table~\ref{tab:sparse-results}, the two algorithms produce essentially the same solution quality: their reconstruction errors, worst sensor residuals, objective values, support recovery F-scores, and stationarity residuals are nearly indistinguishable, and both recover the exact support in $17$ of the $20$ instances. Their average numbers of outer iterations are also very similar. The main difference is computational. Algorithm~\ref{algo:PSG} requires no inner optimization and is approximately one order of magnitude faster than JGS-PPA in terms of mean CPU time. In contrast, JGS-PPA requires, on average, approximately $3.0\times10^4$ inner FISTA iterations per instance to solve the weighted subproblems arising during its scalar searches. Figure~\ref{fig:sparse-residual-time} further illustrates this difference for a representative instance: the two algorithms reach essentially the same final stationarity level, while Algorithm~\ref{algo:PSG} does so in substantially less CPU time. These results indicate that Algorithm~\ref{algo:PSG} attains solution quality comparable to JGS-PPA, while its use of the composite structure of the objective functions leads to substantially simpler and less expensive iterations.

\section{Conclusion}
\label{s:conclu}

We have established a necessary optimality condition and, under additional structural assumptions, sufficient optimality conditions for a broad class of constrained multiobjective optimization problems, where each objective function is the sum of a possibly nonsmooth nonconvex locally Lipschitz function and a differentiable function with a Lipschitz continuous gradient, minus a weakly convex function. Based on these conditions, we have proposed a proximal subgradient algorithm that specifically exploits the structure of the objectives. Theoretical analysis shows that, under mild assumptions, the generated sequence is bounded and all its cluster points are stationary points. Numerical experiments on multi-sensor sparse recovery illustrate that the proposed algorithm attains solution quality comparable to an existing proximal DC algorithm while substantially simplifying the computation of each iteration.

Our work provides a theoretical and computational framework for tackling constrained multiobjective DC-type optimization problems. Future research could explore extensions of this approach to incorporate additional problem-specific structures, such as stochastic or dynamic components, and further investigate its applications to real-world problems.

\paragraph{Acknowledgements.} Part of this work was carried out during the authors' visit to the Vietnam Institute for Advanced Study in Mathematics (VIASM) in 2024. The authors would like to thank VIASM for its financial support and hospitality. The research of NVT is funded by the Hanoi Pedagogical University 2 Foundation for Science and Technology Development under Grant No. HPU2.2025-UT-10. The research of MND is partially supported by the Australian Research Council Discovery Project DP230101749. The authors also thank the anonymous referees for their constructive comments, which helped improve the manuscript.

\end{document}